\documentclass[a4paper, 12pt] {article}
\usepackage[cp1251]{inputenc}
\usepackage[russian, english]{babel}
\usepackage{amsfonts}
\usepackage{mathtext}
\usepackage{cite}

\usepackage{graphicx}

\usepackage{verbatim}
\usepackage{amsmath}
\usepackage{amsthm}
\usepackage{amssymb}
\usepackage{delarray}
\usepackage{mathrsfs}
\usepackage{xcolor}

\begin{document}

\thispagestyle{empty}

\begin{center}
{\Large\bf An inverse Sturm--Liouville-type problem with constant delay and non-zero initial function}

\end{center}

\begin{center}
{\large\bf Sergey Buterin\footnote{Department of Mathematics, Saratov State University, Russia {\it email: buterinsa@sgu.ru}} and Sergey
Vasilev\footnote{Department of Mathematics, Saratov State University, Russia {\it email: i@vasilev-s-v.ru}} }
\end{center}

{\bf Abstract.} We suggest a new statement of the inverse spectral problem for Sturm--Liouville-type operators with constant delay. This
inverse problem consists in recovering the coefficient (often referred to as potential) of the delayed term in the corresponding equation
from the spectra of two boundary value problems with one common boundary condition. However, all studies in this direction focus on the case
of the zero initial function, i.e. they exploit the assumption that the potential vanishes on the corresponding subinterval. In the present
paper, we waive that assumption in favor of a continuously matching initial function, which leads to appearing an additional term with frozen
argument in the equation. For the resulting new inverse problem, we pay a special attention to the situation when one of the spectra is given
only partially. Sufficient conditions and necessary conditions on the corresponding subspectrum for the unique determination of the potential
are obtained, and a constructive procedure for solving the inverse problem is given. In parallel, we obtain the characterization of the
spectra for the zero initial function and the Neumann common boundary condition, which is found to include an additional restriction as
compared with the case of the Dirichlet common condition.


\\

{\large\bf 1. Introduction and main results}
\\

In recent years, there appeared a considerable interest in the inverse problem of recovering an integrable or a square-integrable potential
$q(x)$ in the functional-differential equation
\begin{equation}\label{1}
-y''(x)+q(x)y(x-a)=\lambda y(x),\quad 0<x<\pi,
\end{equation}
with constant delay $a\in(0,\pi)$ given the spectra of two boundary value problems for (\ref{1}) with one common boundary condition (see
\cite{Pik91, FrYur12, Yang, Ign18, BondYur18-1, ButYur19, VPV19, DV19, SatShieh19, WShM19, Dur20, DB21, DB21-2, DB22, BMSh21, WKSh23} and
references therein). For $a=0,$ this problem becomes the classical inverse Sturm--Liouville problem due to Borg \cite{B, novabook, BK19}, but
the nonlocal case $a>0$ requires different approaches. Moreover, it reveals some essentially different effects in solution of the inverse
problem than in the classical situation $a=0.$ Inter alia, as was recently established in \cite{DB21, DB21-2, DB22}, the solution of the
inverse problem may be non-unique if $a\in(0,2\pi/5).$

Various equations with delay have been actively studied from the last century in connection with numerous applications (see, e.g.,
\cite{Mysh, BellCook, Nor, Hal, Skub, AzbMaxRakh, Mur2}). One of specific features of equation (\ref{1}) for $a>0$ is its underdetermination
since the argument of the unknown function $y(x)$ may go beyond the segment $[0,\pi].$ In order to overcome this issue, one can specify an
{\it initial function}, i.e. to assume that $y(x)=f(x)$ for $x\in(-a,0]$ with some known $f(x).$ Alternatively, one can assume that $q(x)=0$
on $(0,a),$ which actually corresponds to specifying $f=0.$ However, we intentionally distinguish these two ways. Indeed, rewriting equation
(\ref{1}) in the form
\begin{equation}\label{nonhom}
-y''(x)+q^+(x)y(x-a)=\lambda y(x)-r(x), \quad 0<x<\pi,
\end{equation}
where $r(x)=q^-(x)f(x-a)$ and
\begin{equation}\label{q+-}
q^-(x)=\left\{\begin{array}{cl}q(x), &x\in(0,a),\\[3mm]
0, & x\in(a,\pi),\end{array}\right. \quad
q^+(x)=\left\{\begin{array}{cl}0, &x\in(0,a), \\[3mm]
q(x), & x\in(a,\pi),\end{array}\right.
\end{equation}
shows that $f\ne0$ leads to a non-homogenous equation, while $f=0$ deals with the corresponding homogenous one. Thus, for posing an
eigenvalue problem, it is natural to choose the latter, i.e. to assume that $q(x)=0$ on $(0,a).$ In particular, the previous studies of
inverse problems for (\ref{1}) were focused namely on this case, i.e. the reconstruction of $q(x)$ was actually carried out only on $(a,\pi)$
since on $(0,a)$ it was a priori assumed to be equal zero.

Meanwhile, admitting a non-zero $f$ also may be appropriate but one should deal with a “linear” initial function, i.e. when $f$ is linearly
dependent on $y$ as, e.g.,
\begin{equation}\label{LinInit}
f(x)=y(0)g(x), \quad -a<x<0.
\end{equation}
This example is quite natural from the point of view of the general theory \cite{Nor} because it ensures a continuous continuation of the
solution $y(x)$ to $[-a,0)$ whenever $g(x)\in C[-a,0]$ and $g(0)=1.$ Such continuation, however, is not always required (see, e.g.,
\cite{AzbMaxRakh}). So one can consider more general forms of the initial function such as, e.g., $f(x)=Ly(x)$ with a linear operator~$L$
acting from $L_2(0,\pi)$ to $L_\infty(-a,0).$ Then for keeping $L$ in frames of a perturbation, a natural requirement would be its relative
compactness \cite{Kato} with respect to the minimal operator of double differentiation. In particular, one can take $Ly(x)=F(y)g(x),$ where
$F(y)$ is a linear functional relatively bounded to that operator. For example, $F(y)=y(b)$ or $F(y)=y'(b)$ for some $b\in[0,\pi].$ We will
focus, however, on the special case~(\ref{LinInit}).

In the present paper, we study the inverse problem for equation (\ref{1}) refusing the usual assumption $q(x)=0$ a.e. on $(0,a)$ in favor of
specifying the “linear” initial function in the form (\ref{LinInit}). Then equation (\ref{1}) can be rewritten with the so-called frozen
argument:
$$
-y''(x)+q^+(x)y(x-a)+p(x)y(0)=\lambda y(x), \quad 0<x<\pi, \quad p(x):=q^-(x)g(x-a).
$$
Since the functions $q^-(x)$ and $g(x-a)$ enter only in their product $p(x),$ they cannot be recovered simultaneously from any spectral
information. Moreover, the reconstruction of $q^-(x)$ on any subinterval $(\alpha,\beta)\subset(0,a)$ can be possible only if $g(x)\ne0$ a.e.
on $(\alpha-a,\beta-a).$ For those reasons, we consider without loss of generality the canonical situation when $g(x)\equiv1.$

For $j=0,1,$ let $\{\lambda_{n,j}\}_{n\ge0}$ be the spectrum of the problem $B_j(q)$ for equation (\ref{1}) with a complex-valued potential
$q(x)\in L_2(0,\pi)$ under the boundary conditions
$$
y'(0)=y^{(j)}(\pi)=0
$$
and under the initial-function condition
$$
y(x)=y(0),\quad -a<x<0.
$$
Consider the following inverse problem.

\medskip
{\bf Inverse Problem 1.} Given $\{\lambda_{n,0}\}_{n\ge0}$ and $\{\lambda_{n,1}\}_{n\ge0},$ find $q(x).$

\medskip
The main results of the present paper (Theorems~1--3) are restricted to the case $a\ge\pi/2.$ In accordance with \cite{DB21-2, DB22}, the
solution of Inverse Problem~1 may be non-unique for $a\in(0,2\pi/5),$ while the case $a\in[2\pi/5,\pi/2)$ requires an additional
investigation. For the future reference, however, we will mark those auxiliary assertions below whose proofs automatically extend to any
wider ranges of $a$ than just $[\pi/2,\pi).$

Everywhere below, one and the same symbol $\{\varkappa_n\}$ will denote {\it different} sequences in $l_2.$ The following theorem gives basic
necessary conditions for the solvability of Inverse Problem~1.

\medskip
{\bf Theorem 1. }{\it For $j=0,1,$ the following asymptotics holds
\begin{equation}\label{ass}
\lambda_{n,j}=\rho_{n,j}^2, \quad \rho_{n,j}=n+\frac{1-j}2 +\frac{\omega}{\pi n}\cos\Big(n+\frac{1-j}2\Big)a +\frac{\varkappa_n}n, \quad
\omega\in{\mathbb C}.
\end{equation}
Here, the constant $\omega$ is determined by the formula
\begin{equation}\label{2.2}
\omega=\frac12\int_a^\pi q^+(x)\,dx.
\end{equation}
Moreover, if the spectra $\{\lambda_{n,0}\}_{n\ge0}$ and $\{\lambda_{n,1}\}_{n\ge0}$ correspond to one and the same $q^-(x),$ then
\begin{equation}\label{char}
i\theta_0(-ir)-\theta_1(-ir)=o(e^{(\pi-a)r}), \quad r\to+\infty,
\end{equation}
where
\begin{equation}\label{A2}
\theta_0(\rho)=\rho(\Delta_0(\rho^2) -\cos\rho\pi) -\omega\sin\rho(\pi-a), \quad \theta_1(\rho)=\Delta_1(\rho^2) +\rho\sin\rho\pi
-\omega\cos\rho(\pi-a),
\end{equation}
while the functions $\Delta_0(\lambda)$ and $\Delta_1(\lambda)$ are determined by the formulae
\begin{equation}\label{A3}
\Delta_0(\lambda)=\prod_{n=0}^\infty\frac{\lambda_{n,0}-\lambda}{(n+1/2)^2}, \quad
\Delta_1(\lambda)=\pi(\lambda_{0,1}-\lambda)\prod_{n=1}^\infty\frac{\lambda_{n,1}-\lambda}{n^2}.
\end{equation} }

Condition (\ref{char}) actually means that Inverse Problem~1 remains overdetermined as in the case $q^-=0$ (see \cite{ButYur19, BMSh21}). As
will be seen below, it is sufficient to specify only one full spectrum and an appropriate part of the other one. For example, we consider
also the following problem.

\medskip
{\bf Inverse Problem 2.} Given $\{\lambda_{n_k,0}\}_{k\in{\mathbb N}}$ and $\{\lambda_{n,1}\}_{n\ge0},$ find $q(x).$

\medskip
Here, $\{n_k\}_{k\in{\mathbb N}}$ is an increasing sequence of non-negative integers. The next theorem gives sufficient conditions as well as
necessary conditions on $\{n_k\}_{k\in{\mathbb N}}$ for the uniqueness of $q(x).$

\medskip
{\bf Theorem 2. }{\it (i) If the system $\sigma_0:=\{\sin(n_k+1/2)x\}_{k\in{\mathbb N}}$ is complete in ${\cal H}:=L_2(0,\pi-a),$ then the
potential $q(x)$ in Inverse Problem~2 is determined uniquely.

(ii) Conversely, if the specification of $\{\lambda_{n_k,0}\}_{k\in{\mathbb N}}$ and $\{\lambda_{n,1}\}_{n\ge0}$ uniquely determines $q(x),$
then the defect of $\sigma_0$ does not exceed $1,$ i.e. $\dim({\cal H}\ominus\sigma_0)\le1.$}

\medskip
Since the system $\{\sin(n+1/2)x\}_{n\ge0}$ is complete in $L_2(0,\pi),$ this theorem, obviously, implies the unique determination of $q(x)$
by both complete spectra as in Inverse Problem~1.

We note that the gap between the sufficient and the necessary conditions in Theorem~2 is actually caused by imposing the common Neumann
boundary condition (at $0).$ By the same reason, the conditions in Theorem~1 do not suffice for the solvability of Inverse Problem~1.

In the case of the Dirichlet common condition, necessary and sufficient conditions for the solvability of the corresponding inverse problem
were obtained in \cite{BMSh21} when $q^-=0.$ Here, we provide such conditions in the same case $q^-=0$ but for the Neumann common condition,
which brings to them an additional item. Specifically, the following theorem holds.

\medskip
{\bf Theorem 3. }{\it Arbitrary complex sequences $\{\lambda_{n,0}\}_{n\ge0}$ and $\{\lambda_{n,1}\}_{n\ge0}$ of the form (\ref{ass}) sharing
one and the same $\omega\in {\mathbb C}$ are the spectra of the problems $B_0(q)$ and $B_1(q),$ respectively, with $q(x)=0$ a.e. on $(0,a)$
if and only if the exponential types of the functions $\theta_0(\rho)$ and $\theta_1(\rho)$ determined by (\ref{A2}) and (\ref{A3}) do not
exceed $\pi-a$ and the following relation is fulfilled:
\begin{equation}\label{A4}
\lambda_{0,1}\prod_{n=1}^\infty\frac{\lambda_{n,1}}{n^2}=\frac{2\omega}\pi.
\end{equation}
}

The latter relation is an additional characterizing condition, which is unnecessary in the Dirichlet case \cite{BMSh21}. We note that the
relevant difference between both cases was pointed out in \cite{DB21} (see Remark~2 therein).

It could be finally mentioned that there recently appeared various studies devoted to the problem of recovering the operator with purely
frozen argument
$$
\ell y:=-y''(x)+q(x)y(b), \quad y^{(\alpha)}(0)=y^{(\beta)}(\pi)=0,
$$
from its spectrum, where $b\in[0,\pi]$ and $\alpha,\beta\in\{0,1\}$ (see \cite{BBV, BV, BK, Wang20, Bond21, MMAS22, DobHry, Kuz22, Kuz23} and
references therein). In particular, its unique solvability depends on the value of $b$ as well as on $\alpha$ and $\beta.$ We note that both
related to Inverse Problem~1 situations: $b=0,$ $\alpha=1,$ $\beta=0$ and  $b=0,$ $\alpha=\beta=1$ belong to the so-called non-generate case,
when the solution is unique (see, e.g., \cite{BBV, BK, MMAS22}).

The paper is organized as follows. In the next section, we construct transformation operators for a fundamental system of solutions of the
homogeneous equation in (\ref{nonhom}), i.e. when $r(x)=0.$ In Section~3, Green's function of the Cauchy problem for the non-homogeneous
equation (\ref{nonhom}) under the zero initial conditions is constructed. In Section~4, we study the characteristics functions of the
problems $B_j(q)$ and prove Theorem~1. Proofs of Theorems~2 and~3 are given in Section~5 along with a constructive procedure for solving the
inverse problems.
\\

{\large\bf 2. Transformation operators}
\\

Let $C(x,\lambda)$ and $S(x,\lambda)$ be solutions of the homogeneous equation in (\ref{nonhom}), i.e. the equation
\begin{equation}\label{hom}
-y''(x)+q^+(x)y(x-a)=\lambda y(x), \quad 0<x<\pi,
\end{equation}
under the initial conditions
$$
C(0,\lambda)=S'(0,\lambda)=1, \quad C'(0,\lambda)=S(0,\lambda)=0.
$$
As in the local case $a=0,$ they form a fundamental system of solutions of equation (\ref{hom}). Throughout the paper, $f'$ and $f^{(j)}$
denote the derivatives with respect to the {\it first} argument:
$$
f'(x_1,\ldots,x_m):=\frac{d}{d x_1}f(x_1,\ldots,x_m), \quad f^{(j)}(x_1,\ldots,x_m):=\frac{d^j}{d x_1^j}f(x_1,\ldots,x_m).
$$

In this section, we obtain representations for the functions $C(x,\lambda)$ and $S(x,\lambda)$ involving the so-called transformation
operators, which connect them with the corresponding solutions of the simplest equation with the zero potential. Specifically, the following
lemma holds.

\medskip
{\bf Lemma 1. }{\it Let $a\ge\pi/2.$ The functions $S(x,\lambda)$ and $C(x,\lambda)$ admit the representations
\begin{equation}\label{trans_S}
S(x,\lambda)=\frac{\sin\rho x}\rho +\int_a^x P(x,t) \frac{\sin\rho(x-t)}\rho\,dt,
\end{equation}
\begin{equation}\label{trans_C}
C(x,\lambda)=\cos\rho x +\int_a^x K(x,t) \cos\rho (x-t)\,dt,
\end{equation}
where $\rho^2=\lambda$ and
\begin{equation}\label{ker S}
P(x,t)=\frac12\int_{\frac{a+t}2}^{x+\frac{a-t}2}q^+(\tau)\,d\tau,
\end{equation}
\begin{equation}\label{ker C}
K(x,t)=\frac12\int_a^{\frac{a+t}2}q^+(\tau)\,d\tau +\frac12\int_a^{x+\frac{a-t}2}q^+(\tau)\,d\tau.
\end{equation}
}

\medskip
{\it Proof.} The assertion for $S(x,\lambda)$ is a particular case of Lemma~1 in \cite{BMSh21}. So we will prove only (\ref{trans_C}) and
(\ref{ker C}). It is easy to see that the Cauchy problem for $C(x,\lambda)$ is equivalent to the integral equation
$$
C(x,\lambda)=\cos\rho x +\int_a^x\frac{\sin\rho(x-t)}\rho q^+(t)C(t-a,\lambda)\,dt.
$$
Taking into account that $a\ge\pi/2,$ we calculate
$$
\int_a^x\frac{\sin\rho(x-t)}\rho q^+(t)C(t-a,\lambda)\,dt =\int_a^x\frac{\sin\rho(x-t)}\rho q^+(t)\cos\rho(t-a)\,dt
$$
$$
=\int_a^x q^+(t)\cos\rho(t-a)\,dt\int_0^{x-t} \cos\rho\tau\,d\tau
$$
$$
=\frac12\int_a^x q^+(t)\,dt\int_0^{x-t} \Big(\cos\rho(t-a+\tau) +\cos\rho(t-a-\tau)\Big)\,d\tau
$$
$$
=\frac12\int_a^x q^+(t)\,dt\int_a^{2(x-t)+a} \cos\rho(x-\tau)\,d\tau =\frac12\int_a^{2x-a} \cos\rho(x-t)\,dt\int_a^{x+\frac{a-t}2}
q^+(\tau)\,d\tau
$$
$$
=\frac12\int_a^x \Big(\int_a^{x+\frac{a-t}2} q^+(\tau)\,d\tau +\int_a^{\frac{a+t}2} q^+(\tau)\,d\tau\Big)\cos\rho(x-t)\,dt,
$$
which finishes the proof. $\hfill\Box$

\medskip
{\bf Remark 1.} While the imposed restriction $a\ge\pi/2$ is vital for (\ref{ker S}) and (\ref{ker C}), representations (\ref{trans_S}) and
(\ref{trans_C}) remain valid also for all smaller $a\ge0$ but with more complicated kernels. In particular, Lemma~1 in \cite{BMSh21} gives an
integral equation for $P(x,t)$ for all $a\in[0,\pi/2).$ Moreover, it extends representation (\ref{trans_S}) to quadratic pencils with two
delays.

\medskip
The following corollary can be easily checked by direct calculations.

\medskip
{\bf Corollary 1. }{\it The following representations hold:
\begin{equation}\label{trans_C-1}
C(x,\lambda)=\cos\rho x +\omega(x)\frac{\sin\rho(x-a)}\rho +\int_a^x K_0(x,t) \frac{\sin\rho(x-t)}\rho\,dt,
\end{equation}
\begin{equation}\label{trans_C'}
C'(x,\lambda)=-\rho\sin\rho x +\omega(x)\cos\rho(x-a) +\int_a^x K_1(x,t) \cos\rho (x-t)\,dt,
\end{equation}
where
\begin{equation}\label{Kernels-1}
\omega(x)=\frac12\int_a^x q^+(t)\,dt, \quad K_j(x,t)=\frac14\Big(q^+\Big(\frac{a+t}2\Big)-(-1)^jq^+\Big(x+\frac{a-t}2\Big)\Big), \;\; j=0,1.
\end{equation}
}
\\

{\large\bf 3. Green's function of the Cauchy operator}
\\

Here, we obtain the solution $z(x,\lambda)=z(x,\lambda;r)$ of the Cauchy problem for the non-homogeneous equation (\ref{nonhom}) with an
arbitrary free term $r(x)$ under the zero initial conditions
\begin{equation}\label{3.0}
z(0,\lambda)=z'(0,\lambda)=0.
\end{equation}
In the next section, we will need representations for $z(\pi,\lambda;q^-)$ and $z'(\pi,\lambda;q^-).$

As in the local case $a=0,$ the function $z(x,\lambda)$ is expected to have the form
\begin{equation}\label{3.1}
z(x,\lambda)=\int_0^x G(x,t,\lambda)r(t)\,dt,
\end{equation}
where $G(x,t,\lambda)$ is the corresponding Green function. Let us find an explicit formula for it.

The following lemma holds for any $a\in[0,\pi].$

\medskip
{\bf Lemma 2. }{\it For each fixed $t\in[0,\pi),$ the function
\begin{equation}\label{3.3-1}
y_t(x):=G(x+t,t,\lambda), \quad 0\le x\le\pi-t,
\end{equation}
is a solution of the Cauchy problem
\begin{equation}\label{3.3}
-y_t''(x)+q_t(x)y_t(x-a)=\lambda y_t(x),\quad 0<x<\pi-t, \quad y_t(0)=0, \quad y_t'(0)=1,
\end{equation}
where
\begin{equation}\label{3.3-2}
q_t(x):=\left\{
\begin{array}{cc}
0,        &   0<x<\min\{a,\pi-t\},\\[3mm]
q^+(x+t), &   a<x<\pi-t.
\end{array}\right.
\end{equation}
}

{\it Proof.} Since the function $G(x,t,\lambda)$ is uniquely determined by the representation (\ref{3.1}), one has the right to impose any
restrictions on it that will finally lead to (\ref{3.1}). In particular, it is natural to assume that $G(x,t,\lambda)$ is sufficiently smooth
and obeys the conditions
\begin{equation}\label{3.3-3}
G(x,x,\lambda)=0, \quad G'(x,x,\lambda)=1.
\end{equation}
Then substituting (\ref{3.1}) into (\ref{nonhom}) and taking the arbitrariness of $r(x)$ into account, we obtain the relations
$$
\begin{array}{rl}
-G''(x,t,\lambda) =\lambda G(x,t,\lambda), & 0<t<x<a,\\[3mm]
-G''(x,t,\lambda) +q^+(x)G(x-a,t,\lambda) =\lambda G(x,t,\lambda), & 0<t<x-a<\pi-a,\\[3mm]
-G''(x,t,\lambda) =\lambda G(x,t,\lambda), & 0<x-a<t<x<\pi,
\end{array}
$$
which along with (\ref{3.3-3}), in turn, guaranty that (\ref{3.1}) is a solution of the problem (\ref{nonhom}) and (\ref{3.0}).

Substituting $x+t$ into the above three relations instead of $x,$ we get
\begin{equation}\label{3.2}
-G''(x+t,t,\lambda) =\lambda G(x+t,t,\lambda), \quad 0<x<a-t<a,
\end{equation}
\begin{equation}\label{3.2-1}
-G''(x+t,t,\lambda) +q^+(x+t)G(x+t-a,t,\lambda) =\lambda G(x+t,t,\lambda), \quad a<x<\pi-t<\pi,
\end{equation}
\begin{equation}\label{3.2-2}
-G''(x+t,t,\lambda) =\lambda G(x+t,t,\lambda), \quad \max\{0,a-t\}<x<\min\{a,\pi-t\}.
\end{equation}
Combining (\ref{3.2}) and (\ref{3.2-2}) and taking (\ref{3.3-1}) into account we rewrite:
$$
-y_t''(x)=\lambda y_t(x), \quad 0<x<\min\{a,\pi-t\},
$$
while (\ref{3.2-1}) takes the form
$$
-y_t''(x)+q^+(x+t)y_t(x-a)=\lambda y_t(x), \quad a<x<\pi-t<\pi.
$$
Using the designation (\ref{3.3-2}) along with initial conditions (\ref{3.3-3}), we arrive at (\ref{3.3}).

Finally, note that, after solving the Cauchy problem (\ref{3.3}) in the standard way (see, e.g., \cite{DB22}), it is easy to see that
$G(x,t,\lambda)$ is a continuous function with respect to all arguments. Hence, the integral in (\ref{3.1}) exists and gives a solution to
the Cauchy problem (\ref{nonhom}) and (\ref{3.0}). $\hfill\Box$

\medskip
{\bf Lemma 3. }{\it Let $a\ge\pi/2.$ Then the following representations hold:
\begin{equation}\label{3.4}
G(x,t,\lambda)=\frac{\sin\rho(x-t)}\rho, \quad \max\{0,x-a\}\le t\le x\le \pi,
\end{equation}
and
\begin{equation}\label{3.5}
G(x,t,\lambda)=\frac{\sin\rho(x-t)}\rho +\frac12\int_{a+t}^x \frac{\sin\rho(x-\tau)}\rho\,d\tau \int_{\frac{a+t+\tau}2}^{x+\frac{a+t-\tau}2}
q^+(\eta)\,d\eta
\end{equation}
whenever $0\le t\le x-a\le \pi-a.$ }

\medskip
{\it Proof.} By virtue of (\ref{3.3}) and Lemma~1, we have the representation
$$
y_t(x)=\frac{\sin\rho x}\rho +\frac12\int_a^x \frac{\sin\rho(x-\tau)}\rho\,d\tau \int_{\frac{a+\tau}2}^{x+\frac{a-\tau}2} q_t(\eta)\,d\eta,
\quad 0\le x\le\pi-t,
$$
which, in accordance with (\ref{3.3-1}) and (\ref{3.3-2}), leads to (\ref{3.4}) and (\ref{3.5}). $\hfill\Box$

\medskip
By substituting (\ref{3.4}) and (\ref{3.5}) into (\ref{3.1}) and changing the order of integration, we obtain
\begin{equation}\label{3.9}
z(x,\lambda)=\int_0^x \Big(r(t) +\frac12\int_0^{t-a}r(\tau)\,d\tau \int_{\frac{a+t+\tau}2}^{x+\frac{a+\tau-t}2} q^+(\eta)\,d\eta\Big)
\frac{\sin\rho(x-t)}\rho\,dt, \quad 0\le x\le\pi,
\end{equation}
where $r(x)=0$ for $x<0.$

Further, differentiating (\ref{3.4}) and (\ref{3.5}) with respect to $x,$ we arrive at the formulae
$$
G'(x,t,\lambda)=\cos\rho(x-t), \quad \max\{0,x-a\}\le t\le x\le \pi,
$$
and
$$
G'(x,t,\lambda)=\cos\rho(x-t) +\frac12\int_{a+t}^x \Big(\int_{\frac{a+t+\tau}2}^x q^+(\eta)\,d\eta +\int_{x+\frac{a+t-\tau}2}^x
q^+(\eta)\,d\eta\Big)\cos\rho(x-\tau)\,d\tau
$$
as soon as $0\le t\le x-a\le \pi-a.$ Substituting them into
$$
z'(x,\lambda)=\int_0^x G'(x,t,\lambda)r(t)\,dt,
$$
we analogously obtain the representation
\begin{equation}\label{3.10}
\begin{array}{c} \displaystyle z'(x,\lambda)=\int_0^x \Big(r(t) + \frac12\int_0^{t-a}\Big(\int_{\frac{a+t+\tau}2}^x q^+(\eta)\,d\eta
\qquad\qquad\qquad\qquad\qquad\qquad\qquad\qquad\;\;
\\[5mm]
\displaystyle
 \;\;\qquad\qquad\qquad\qquad\qquad\qquad\qquad\qquad +\int_{x+\frac{a+\tau-t}2}^x q^+(\eta)\,d\eta\Big)r(\tau)\,d\tau\Big)\cos\rho(x-t)\,dt.
\end{array}
\end{equation}
\\

{\large\bf 4. Characteristic functions}
\\

Consider the entire functions
\begin{equation}\label{4.1}
\Delta_j(\lambda):=C^{(j)}(\pi,\lambda) +z^{(j)}(\pi,\lambda;q^-), \quad j=0,1.
\end{equation}
The next lemma holds for any $a\in[0,\pi].$

\medskip
{\bf Lemma 4. }{\it For $j=0,1,$ eigenvalues of the problem $B_j(q)$ coincide with zeros of $\Delta_j(\lambda).$}

\medskip
{\it Proof.} Since the sum $C(x,\lambda) +z(x,\lambda;q^-)$ cannot be identically zero, any zero of $\Delta_j(\lambda)$ is an eigenvalue of
the problem $B_j(q),$ which under our settings, in turn, has the form
\begin{equation}\label{4.2}
-y''(x)+q^+(x)y(x-a)+q^-(x)y(0)=\lambda y(x), \quad y'(0)=y^{(j)}(\pi)=0.
\end{equation}
Conversely, let $\lambda$ be an eigenvalue of $B_j(q),$ and let $y(x)$ be the corresponding eigenfunction, i.e. a nontrivial solution of
(\ref{4.2}). Then $y(0)\ne0$ since, obviously, $y(x)\equiv0$ otherwise. Without loss of generality, one can assume that $y(0)=1,$ which will
imply $y(x)=C(x,\lambda) +z(x,\lambda;q^-)$ due to uniqueness of solution of the Cauchy problem. Hence, $\Delta_j(\lambda)=y^{(j)}(\pi)=0.$
$\hfill\Box$

\medskip
As usual, we call $\Delta_j(\lambda)$ {\it characteristic function} of the problem $B_j(q).$ The following lemma based on the two preceding
sections gives representations for both characteristic functions.

\medskip
{\bf Lemma 5. }{\it The characteristic functions admit the representations
\begin{equation}\label{CharInt-1}
\Delta_0(\lambda)=\cos\rho\pi +\omega\frac{\sin\rho(\pi-a)}\rho +\int_0^\pi w_0(x)\frac{\sin\rho x}\rho\,dx, \quad w_0(x)\in L_2(0,\pi),
\end{equation}
\begin{equation}\label{CharInt-2}
\Delta_1(\lambda)=-\rho\sin\rho\pi +\omega\cos\rho(\pi-a) +\int_0^\pi w_1(x)\cos\rho x\,dx, \quad w_1(x)\in L_2(0,\pi).
\end{equation}
Moreover, the constant $\omega$ is determined by (\ref{2.2}) and
\begin{equation}\label{0-a}
w_0(\pi-x)=w_1(\pi-x)=q^-(x), \quad 0<x<a,
\end{equation}
while for $a<x<\pi:$
\begin{equation}\label{a-pi-0}
w_0(\pi-x)=\frac14\Big(q^+\Big(\frac{a+x}2\Big)-q^+\Big(\pi+\frac{a-x}2\Big)\Big)+\frac12 \int_0^{x-a}
q^-(t)\,dt\int_{\frac{a+x+t}2}^{\pi+\frac{a+t-x}2}q^+(\tau)\,d\tau,
\end{equation}
$$
w_1(\pi-x) =\frac14\Big(q^+\Big(\frac{a+x}2\Big) +q^+\Big(\pi+\frac{a-x}2\Big)\Big) \qquad\qquad\qquad\qquad\qquad\qquad\qquad\qquad\qquad
$$
\begin{equation}\label{a-pi-1}
\quad\qquad\qquad\qquad\qquad\qquad\qquad +\frac12 \int_0^{x-a}\Big( \int_{\frac{a+x+t}2}^\pi q^+(\tau)\,d\tau +\int_{\pi+\frac{a+t-x}2}^\pi
q^+(\tau)\,d\tau\Big) q^-(t)\,dt.
\end{equation}
}

\medskip
{\it Proof.} Substituting $x=\pi$ into (\ref{trans_C-1}) and (\ref{trans_C'}) and using (\ref{2.2}) and (\ref{Kernels-1}), we obtain
\begin{equation}\label{C-pi}
C(\pi,\lambda)=\cos\rho\pi +\omega\frac{\sin\rho(\pi-a)}\rho +\int_0^{\pi-a} u_0(x) \frac{\sin\rho x}\rho\,dx,
\end{equation}
\begin{equation}\label{C'-pi}
C'(\pi,\lambda)=-\rho\sin\rho\pi +\omega\cos\rho(\pi-a) +\int_0^{\pi-a} u_1(x) \cos\rho x\,dx,
\end{equation}
where
\begin{equation}\label{u_0}
u_j(\pi-x)=K_j(\pi,x) =\frac14\Big(q^+\Big(\frac{a+x}2\Big)-(-1)^jq^+\Big(\pi+\frac{a-x}2\Big)\Big), \;\; a<x<\pi, \;\; j=0,1.
\end{equation}
Further, substituting $r=q^-$ and $x=\pi$ into (\ref{3.9}) and (\ref{3.10}), we arrive at
\begin{equation}\label{z-pi}
z(\pi,\lambda;q^-)=\int_0^\pi v_0(x) \frac{\sin\rho x}\rho\,dx, \quad z'(\pi,\lambda;q^-)=\int_0^\pi v_1(x) \cos\rho x\,dx,
\end{equation}
where
\begin{equation}\label{z_0-a}
v_0(\pi-x)=v_1(\pi-x)=q^-(x), \quad 0<x<a,
\end{equation}
\begin{equation}\label{z_a-pi-0}
v_0(\pi-x)=\frac12 \int_0^{x-a} q^-(t)\,dt\int_{\frac{a+x+t}2}^{\pi+\frac{a+t-x}2}q^+(\tau)\,d\tau, \quad a<x<\pi,
\end{equation}
\begin{equation}\label{z_a-pi-1}
v_1(\pi-x) =\frac12 \int_0^{x-a}\Big( \int_{\frac{a+x+t}2}^\pi q^+(\tau)\,d\tau +\int_{\pi+\frac{a+t-x}2}^\pi q^+(\tau)\,d\tau\Big)
q^-(t)\,dt, \quad a<x<\pi.
\end{equation}
According to (\ref{4.1}), (\ref{C-pi}), (\ref{C'-pi}) and (\ref{z-pi}), we get (\ref{CharInt-1}) and (\ref{CharInt-2}) with
\begin{equation}\label{4.3}
w_j(x)=u_j(x)+v_j(x), \quad j=0,1,
\end{equation}
where $u_0(x)=u_1(x)=0$ on $(\pi-a,\pi).$ Finally, substituting (\ref{u_0}) and (\ref{z_0-a})--(\ref{z_a-pi-1}) into (\ref{4.3}), we arrive
at (\ref{0-a})--(\ref{a-pi-1}). $\hfill\Box$

\medskip
In the rest part of this section, we provide auxiliary facts about arbitrary functions of the form (\ref{CharInt-1}) and (\ref{CharInt-2}),
and give the proof of Theorem~1.

Lemmas~6--8 below are valid for any fixed $a\in[0,2\pi].$ By the standard approach (see, e.g., \cite{novabook, But22}) involving Rouch\'e's
theorem, one can prove the following assertion.

\medskip
{\bf Lemma 6. }{\it For $j=0,1,$ any $\Delta_j(\lambda)$ has infinitely many zeros $\{\lambda_{n,j}\}_{n\ge0}$ of the form (\ref{ass}).}

\medskip
The next assertion for $a=0$ can be found in \cite{novabook} but the proof does not depend on the value of $a$ as soon as it ranges within
$[0,2\pi].$

\medskip
{\bf Lemma 7. }{\it Any functions of the forms (\ref{CharInt-1}) and (\ref{CharInt-2}) are determined by their zeros uniquely. Moreover, the
representations in (\ref{A3}) hold.}

\medskip
Now, we are in position to give the proof of Theorem~1.

\medskip
{\bf Proof of Theorem~1.} The asymptotics (\ref{ass}) is a direct corollary of Lemmas~5 and~6. It remains to note that, by virtue of
(\ref{A2}), (\ref{CharInt-1}) and (\ref{CharInt-2}) along with Lemma~7, we have
$$
i\theta_0(\rho)-\theta_1(\rho)=i\int_0^\pi w_0(x)\sin\rho x\,dx -\int_0^\pi w_1(x)\cos\rho x\,dx=\frac{\theta_+(\rho)-\theta_-(\rho)}2,
$$
where, according to (\ref{0-a}),
$$
\theta_+(\rho)=\int_0^{\pi-a} (w_0-w_1)(x)\exp(i\rho x)\,dx, \quad \theta_-(\rho)=\int_0^\pi (w_0+w_1)(x)\exp(-i\rho x)\,dx,
$$
which implies (\ref{char}). $\hfill\Box$

\medskip
Statements analogous to the next lemma are often used for finding necessary and sufficient conditions for the solvability of inverse
problems, i.e. a characterization of the spectral data (see Remark~2 in \cite{But22}). For its proof, we will follow a new simple idea
suggested in~\cite{But22}.

\medskip
{\bf Lemma 8. }{\it For $j=0,1,$ let $\{\lambda_{n,j}\}_{n\ge0}$ be arbitrary complex sequences of the form~(\ref{ass}). Then the function
$\Delta_j(\lambda)$ constructed by the corresponding formula in (\ref{A3}) has the form (\ref{CharInt-1}) or (\ref{CharInt-2}),
respectively.}

\medskip
{\it Proof.} Since the assertion of the lemma for $j=0$ formally follows from Lemma~6 in \cite{BMSh21}, we focus on the case $j=1.$ Let a
sequence $\{\lambda_{n,1}\}_{n\ge0}$ of the form~(\ref{ass}) be given. First, let all values $\lambda_{n,1}$ be distinct and
$\lambda_{0,1}=0.$ Denote $\rho_{-n,1}:=-\rho_{n,1}$ for $n\ge1.$ By virtue of Lemma~2 in \cite{But22}, the system
$\{\exp(i\rho_{n,1}x)\}_{n\in{\mathbb Z}}$ is a Riesz basis in $L_2(-\pi,\pi).$ Moreover, the asymptotics (\ref{ass}) implies
$\{\theta(\rho_{n,1})\}_{n\in{\mathbb Z}}\in l_2,$ where $\theta(\rho):=\rho\sin\rho\pi -\omega\cos\rho(\pi-a)$ and $\omega$ is as
in~(\ref{ass}). Hence, there exists a unique function $W_1(x)\in L_2(-\pi,\pi)$ obeying the relations
$$
\theta(\rho_{n,1})=\int_{-\pi}^\pi W_1(x)\exp(i\rho_{n,1}x)\,dx, \quad n\in{\mathbb Z}.
$$
Obviously,  $W_1(x)$ is even. Thus, $\lambda_{n,1}=(\rho_{n,1})^2,$ $n\ge0,$ are zeros of the function $\Delta_1(\lambda)$ determined by
(\ref{CharInt-2}) with $w_1(x)=2W_1(x).$ By Lemma~6, $\Delta_1(\lambda)$ has no other zeros, while by Lemma~7, it admits the second
representation in (\ref{A3}), which finishes the proof for a simple sequence $\{\lambda_{n,1}\}_{n\ge0}$ containing a zero element.

For the general case, it is sufficient to note that multiplying $\Delta_1(\lambda)$ with any function
$$
h(\lambda):=\prod_{n\in A} \frac{\lambda-\tilde\lambda_{n,1}}{\lambda-\lambda_{n,1}}, \quad A\subset{\mathbb N}\cup\{0\}, \quad \#A<\infty,
$$
preserves the form (\ref{CharInt-2}) and changes only $w_1(x).$ Indeed, we have
$$
h(\lambda)\Delta_1(\lambda)=-\rho\sin\rho\pi +\omega\cos\rho(\pi-a) +H(\lambda),
$$
where
$$
H(\lambda)= (1-h(\lambda))\Big(\rho\sin\rho\pi -\omega\cos\rho(\pi-a)\Big) +h(\lambda)\int_0^\pi w_1(x)\cos\rho x\,dx.
$$
The function $H(\lambda)$ is entire as soon as $\lambda_{n,1}$ are zeros of $\Delta_1(\lambda).$ Moreover, in the $\rho$-plane, we,
obviously, have $H(\rho^2)\in L_2(-\infty,+\infty)$ and $H(\rho^2)=o(\exp(|{\rm Im}\,\rho|\pi))$ as $\rho\to\infty.$ Thus, by virtue of the
Paley--Wiener theorem (see, e.g., \cite{Lev}), it has the form
$$
H(\lambda)= \int_0^\pi \tilde w_1(x)\cos\rho x\,dx, \quad \tilde w_1(x)\in L_2(0,\pi),
$$
which finishes the proof completely. $\hfill\Box$

\medskip
Finally, let us give one more auxiliary assertion, which will be used in the proof of Theorem~2. Let $\{n_k\}_{k\in{\mathbb N}}$ be an
increasing sequence of non-negative integers. Without loss of generality, assume that multiple elements in the subspectrum
$\{\lambda_{n_k,0}\}_{k\in{\mathbb N}}$ are neighboring, i.e.
$$
\lambda_{n_k,0}=\lambda_{n_{k+1},0}=\ldots=\lambda_{n_{k+m_k-1},0},
$$
where $m_k$ is the multiplicity of the value $\lambda_{n_k,0}$ in this subspectrum.  Put
$$
{\cal S}:=\{1\}\cup\{k:\lambda_{n_k,0}\ne\lambda_{n_{k-1},0},\,k\ge 2\}
$$
and consider the functional system $\sigma:=\{s_n(x)\}_{n\in{\mathbb N}},$ where
$$
s_{k+\nu}(x)=(n_k+1)\frac{d^\nu}{d\lambda^\nu}\frac{\sin\rho x}{\rho}\Big|_{\lambda=\lambda_{n_k,0}}, \quad k\in {\cal S}, \quad
\nu=\overline{0,m_k-1}.
$$

\medskip
{\bf Lemma 9. }{\it The system $\sigma$ is complete (is a Riesz basis) in ${\cal H}_b:=L_2(0,b)$ if and only if so is the system $\sigma_0.$
Moreover, they have equal defects, i.e. $\dim({\cal H}_b\ominus\sigma_0)=\dim({\cal H}_b\ominus\sigma).$}

\medskip
{\it Proof.} The first assertion of the lemma coincides with the second assertion of Lemma~1 in \cite{ButYur19}. For proving the second one,
let there exist $d$ linearly independent entire functions $h_\nu(\lambda),$ $\nu=\overline{1,d},$ of the form
$$
h_\nu(\lambda)=\int_0^bf_\nu(x)\frac{\sin\rho x}\rho\,dx, \quad f_\nu(x)\in L_2(0,b),
$$
whose zeros have the common part $\{(n_k+1/2)^2\}_{k\in{\mathbb N}}.$ Consider the meromorphic function
$$
F(\lambda):=\prod_{k=1}^\infty\frac{\lambda_{n_k,0}-\lambda}{(n_k+1/2)^2-\lambda}.
$$
Then the function $\tilde h_\nu(\lambda):=F(\lambda)h_\nu(\lambda)$ also has the form
\begin{equation}\label{4.8}
\tilde h_\nu(\lambda)=\int_0^b\tilde f_\nu(x)\frac{\sin\rho x}\rho\,dx, \quad \tilde f_\nu(x)\in L_2(0,b).
\end{equation}
Indeed, as in the proof of Lemma~2 in \cite{But22}, one can show that $|F(\rho^2)|<C_\delta$ whenever
$$
|\rho\pm (n_k+1/2)|\ge\delta, \quad k\in{\mathbb N},
$$
for each fixed $\delta>0.$ Obviously, the function $\tilde h_\nu(\lambda),$ after removing the singularities, is entire and, by the latter
estimate, we have $|\rho\tilde h_\nu(\rho^2)|\le C_\delta|\rho h_\nu(\rho^2)|\in L_2(-\infty+i\delta,\infty+i\delta)$ in the $\rho$-plane.
Moreover, the maximum modulus principle for analytic functions implies $\rho\tilde h_\nu(\rho^2)=o(\exp(|{\rm Im}\,\rho|b))$ as
$\rho\to\infty.$ Hence, by the Paley--Wiener theorem \cite{Lev}, we have (\ref{4.8}).

The constructed functions $\tilde h_\nu(\lambda),$ $\nu=\overline{1,d},$ are, obviously, linearly independent and their zeros have the common
part $\{\lambda_{n_k,0}\}_{k\in{\mathbb N}}$ with account of multiplicity. Hence, $\dim({\cal H}\ominus\sigma_0)\le\dim({\cal
H}\ominus\sigma).$ The opposite inequality can be proved similarly. $\hfill\Box$
\\

{\large\bf 5. Solution of the inverse problems}
\\

When the functions $w_0(x)$ and $w_1(x)$ are specified, relations (\ref{0-a})--(\ref{a-pi-1}) can be considered as a nonlinear integral
equation with respect to $q(x)=q^-(x)+q^+(x).$ The following lemma actually implies its unique solvability.

\medskip
{\bf Lemma 10. }{\it For any functions $w_0(x),w_1(x),q^-(x)\in L_2(0,\pi-a),$ the linear system consisting of (\ref{a-pi-0}) and
(\ref{a-pi-1}) has a unique solution  $q^+(x)\in L_2(a,\pi).$}

\medskip
{\it Proof.} Adding up equations (\ref{a-pi-0}) and (\ref{a-pi-1}) and then subtracting one from the other, we get
$$
\left.\begin{array}{l}
\displaystyle 2(w_1+w_0)(\pi-x) =q^+\Big(\frac{a+x}2\Big) +2\int_0^{x-a} q^-(t)\,dt \int_{\frac{a+x+t}2}^\pi q^+(\tau)\,d\tau,\\[5mm]
\displaystyle 2(w_1-w_0)(\pi-x) =q^+\Big(\pi+\frac{a-x}2\Big) +2\int_0^{x-a} q^-(t)\,dt \int_{\pi+\frac{a+t-x}2}^\pi q^+(\tau)\,d\tau,
\end{array}\right\} \quad a<x<\pi.
$$
Changing the variable, we arrive at the relations
$$
2(w_1+w_0)(\pi+a-2x) =q^+(x) +2\int_0^{2(x-a)} q^-(t)\,dt \int_{x+\frac{t}2}^\pi q^+(\tau)\,d\tau, \quad a<x<\frac{a+\pi}2,
$$
$$
2(w_1-w_0)(2x-\pi-a) =q^+(x) +2\int_0^{2(\pi-x)} q^-(t)\,dt \int_{x+\frac{t}2}^\pi q^+(\tau)\,d\tau,  \quad \frac{a+\pi}2<x<\pi.
$$
Then changing the order of integration in the last two formulae we obtain the system
$$
2(w_1+w_0)(\pi+a-2x) =q^+(x) +2\int_x^{2x-a} q^+(t)\,dt \int_0^{2(t-x)} q^-(\tau)\,d\tau \qquad\qquad\qquad\quad
$$
$$
\qquad\qquad\qquad\qquad\qquad +2\int_{2x-a}^\pi q^+(t)\,dt \int_0^{2(x-a)} q^-(\tau)\,d\tau, \quad a<x<\frac{a+\pi}2,
$$
$$
2(w_1-w_0)(2x-\pi-a) =q^+(x) +2\int_x^\pi q^+(t)\,dt \int_0^{2(t-x)} q^-(\tau)\,d\tau,  \quad \frac{a+\pi}2<x<\pi.
$$
Using the designations
\begin{equation}\label{5.4}
W(x):=\left\{\begin{array}{c}
\displaystyle 2(w_1+w_0)(\pi+a-2x), \quad a<x<\frac{a+\pi}2,\\[3mm]
\displaystyle 2(w_1-w_0)(2x-\pi-a),  \quad \frac{a+\pi}2<x<\pi,
\end{array}\right.
\end{equation}
\begin{equation}\label{5.5}
Q(x,t):=\left\{\begin{array}{r}
\displaystyle 2\int_0^{2(x-a)} q^-(\tau)\,dt, \quad a<2x-a<t<\pi,\\[3mm]
\displaystyle 2\int_0^{2(t-x)} q^-(\tau)\,dt, \quad a<x<t<\min\{2x-a,\pi\},
\end{array}\right.
\end{equation}
one can rewrite the latter system as the Volterra integral equation
\begin{equation}\label{5.6}
W(x)=q^+(x) +\int_x^\pi Q(x,t)q^+(t)\,dt,  \quad a<x<\pi,
\end{equation}
which possesses a unique solution $q^+(x)\in L_2(a,\pi).$ $\hfill\Box$

\medskip
{\bf Proof of Theorem~2.} First of all, note that, due to (\ref{ass}), the value $\omega$ is always determined by specifying
$\{\lambda_{n,1}\}_{n\ge0}$ via the formula
\begin{equation}\label{5.7}
\omega=\pi\lim_{k\to\infty} \tilde n_k\frac{\rho_{\tilde n_k,1}-\tilde n_k}{\cos \tilde n_ka},
\end{equation}
where the natural sequence $\{\tilde n_k\}$ is chosen so that $|\cos \tilde n_ka|\ge c>0.$ Alternatively, in accordance with
(\ref{CharInt-2}), one can use the formula
\begin{equation}\label{5.8}
\omega=\lim_{n\to\infty} \Big(\Delta_1(\xi_n^2) +\xi_n\sin\xi_n\pi\Big), \quad \xi_n=\frac{2\pi n}{\pi-a},
\end{equation}
where $\Delta_1(\lambda)$ is constructed by the second representation in (\ref{A3}).

(i) Let the system $\sigma_0$ be complete in ${\cal H.}$ Since, according to Lemma~7, the characteristic function $\Delta_1(\lambda)$ is
uniquely determined by its zeros, so is also $w_1(x)$ in (\ref{CharInt-2}). By virtue of~(\ref{0-a}), the function $w_0(x)$ coincides with
$w_1(x)$ a.e. on $(\pi-a,\pi),$ i.e. it becomes known too.

By differentiating (\ref{CharInt-1}) $\nu=\overline{0,m_k-1}$ times and substituting $\lambda=\lambda_{n_k,0}$ for $k\in{\cal S},$ we arrive
at the relations
\begin{equation}\label{5.9}
\beta_n=\int_0^{\pi-a} w_0(x)s_n(x)\,dx, \quad n\in{\mathbb N},
\end{equation}
where $m_k,$ ${\cal S}$ and $s_n(x)$ were defined before Lemma~9 and
\begin{equation}\label{5.10}
\beta_{k+\nu}=-(n_k+1)\left.\frac{d^\nu}{d\lambda^\nu}\Big(\cos\rho\pi
+\omega\frac{\sin\rho(\pi-a)}\rho+\gamma(\lambda)\Big)\right|_{\lambda=\lambda_{n_k,0}}, \;\; k\in{\cal S}, \;\; \nu=\overline{0,m_k-1},
\end{equation}
\begin{equation}\label{5.11}
\gamma(\lambda)=\int_{\pi-a}^\pi w_1(x)\frac{\sin\rho x}\rho\,dx.
\end{equation}
Hence, by virtue of Lemma~9, the function $w_0(x)$ is determined uniquely also on  $(0,\pi-a).$ Thus, it remains to recall representations
(\ref{q+-}) and (\ref{0-a}) as well as to apply Lemma~10.

(ii) Assume that $q(x)$ is uniquely determined by $\{\lambda_{n_k,0}\}_{k\in{\mathbb N}}$ and $\{\lambda_{n,1}\}_{n\ge0}$ and, to the
contrary, that $\dim({\cal H}\ominus\sigma_0)>1.$ Then, according to Lemma~9, we have $\dim({\cal H}\ominus\sigma)>1,$ i.e. there exist at
least two linearly independent functions $f_1(x),f_2(x)\in L_2(0,\pi-a)$ such that
\begin{equation}\label{5.12}
\int_0^{\pi-a}f_\nu(x)s_n(x)\,dx =0, \quad n\in{\mathbb N}, \quad \nu=1,2.
\end{equation}
Consider a function
\begin{equation}\label{5.12-1}
\tilde w_0(x)=w_0(x)+\alpha_1f_1(x)+\alpha_2f_2(x)\in L_2(0,\pi-a), \quad \alpha_1,\alpha_2\in{\mathbb C}.
\end{equation}

Let $\tilde q^+(x)$ become a solution of equation (\ref{5.6}) after replacing its left-hand side with
\begin{equation}\label{5.13}
\tilde W(x)=W(x)+\alpha_1F_1(x)+\alpha_2F_2(x),
\end{equation}
where
\begin{equation}\label{5.14}
F_\nu(x):=\left\{\begin{array}{r}
\displaystyle 2f_\nu(\pi+a-2x), \quad a<x<\frac{a+\pi}2,\\[3mm]
\displaystyle -2f_\nu(2x-\pi-a),  \quad \frac{a+\pi}2<x<\pi.
\end{array}\right.
\end{equation}
In other words, $\tilde q^+(x)=q^+(x)+\alpha_1g_1(x)+\alpha_2g_2(x),$ where
$$
g_\nu(x)=F_\nu(x) +\int_x^\pi \tilde Q(x,t)F_\nu(t)\,dt,  \quad \nu=1,2,
$$
while $\tilde Q(x,t)$ is the resolvent kernel for the kernel $Q(x,t).$ Choose $\alpha_1$ and $\alpha_2$ so that they would not vanish
simultaneously and
\begin{equation}\label{5.15}
\frac12\int_a^\pi \tilde q^+(x)\,dx=\omega.
\end{equation}
Since the functions $F_1(x)$ and $F_2(x)$ are linearly independent, so are $g_1(x)$ and $g_2(x).$ Hence, $\tilde q^+\ne q^+.$ Continue
$\tilde q^+(x)$ to $(0,a)$ as zero and consider the function $\tilde q(x)=q^-(x)+\tilde q^+(x).$ By virtue of (\ref{5.15}) and Lemma~5, the
characteristic functions $\tilde\Delta_0(\lambda)$ and $\tilde\Delta_1(\lambda)$ of the problems $B_0(\tilde q)$ and $B_1(\tilde q),$
respectively, have the forms
$$
\tilde\Delta_0(\lambda)=\cos\rho\pi +\omega\frac{\sin\rho(\pi-a)}\rho +\int_0^\pi \tilde w_0(x)\frac{\sin\rho x}\rho\,dx,
$$
$$
\tilde\Delta_1(\lambda)=-\rho\sin\rho\pi +\omega\cos\rho(\pi-a) +\int_0^\pi \tilde w_1(x)\cos\rho x\,dx.
$$
Comparing (\ref{5.4}), (\ref{5.6}), (\ref{5.13}) and (\ref{5.14}), one can see that $\tilde w_1(x)=w_1(x).$ Hence, the spectra of $B_1(\tilde
q)$ and $B_1(q)$ coincide. Moreover, according to (\ref{5.9})--(\ref{5.12-1}), the sequence $\{\lambda_{n_k,0}\}_{k\in{\mathbb N}}$ is a
subsequence of zeros of $\tilde\Delta_0(\lambda).$ Hence, this sequence is a subspectrum of the problem $B_0(\tilde q).$ Thus, we obtained
another potential $\tilde q\ne q$ with the same spectral data $\{\lambda_{n_k,0}\}_{k\in{\mathbb N}}$ and $\{\lambda_{n,1}\}_{n\ge0}$ as $q$
has. This contradiction finishes the proof. $\hfill\Box$

\medskip
Now, we are in position to give a constructive procedure for solving Inverse Problem~1.

\medskip
{\bf Algorithm 1. }{\it Let the spectra $\{\lambda_{n,0}\}_{n\ge0}$ and $\{\lambda_{n,1}\}_{n\ge0}$ be given. Then:

(i) Construct the functions $\Delta_0(\lambda)$ and $\Delta_1(\lambda)$ by the formulae in (\ref{A3});

(ii) Find the value $\omega$ by (\ref{5.7}) or (\ref{5.8});

(iii) Calculate the functions $w_0(x)$ and $w_1(x)$ in (\ref{CharInt-1}) and (\ref{CharInt-2}) by inverting the corresponding Fourier
transforms:
$$
w_0(x)=\frac2\pi\sum_{n=1}^\infty a_n \sin nx, \quad w_1(x)=\frac2\pi\sum_{n=0}^\infty b_n \cos nx,
$$
where
$$
a_n=n(\Delta_0(n^2)-(-1)^n)+\omega(-1)^n \sin na, \;\; n\ge1, \quad b_n=\Delta_1(n^2)-\omega(-1)^n \cos na, \;\; n\ge0;
$$

(iv) Find $q^-(x)\in L_2(0,a)$ by any relation in (\ref{0-a}) and put $q^-(x)=0$ for $x\in(a,\pi);$

(v) Construct the functions $W(x)$ and $Q(x,t)$ by the formulae (\ref{5.4}) and (\ref{5.5}), respectively, and find $q^+(x)\in L_2(a,\pi)$ by
solving the Volterra integral equation (\ref{5.6});

(vi) Finally, construct $q(x)=q^-(x)+q^+(x),$ where $q^+(x)=0$ on $(0,a).$}

\medskip
This algorithm can be easily extended to Inverse Problem~2 if $\{\sin(n_k+1/2)x\}_{k\in{\mathbb N}}$ is a Riesz basis in $L_2(0,\pi-a).$
Then, by virtue of Lemma~9, so is the system $\{s_n(x)\}_{n\in{\mathbb N}}.$ Therefore, on step (iii), the function $\omega_0(x)$ can be
constructed in accordance with (\ref{5.9}) by  the formula
$$
w_0(x)=\sum_{n=1}^\infty \beta_ns_n^*(x) ,\quad 0<x<\pi-a,
$$
where the coefficients $\beta_n$ are determined by relations (\ref{5.10}) and (\ref{5.11}), while $\{s_n^*(x)\}_{n\in{\mathbb N}}$ is the
biorthogonal basis to the basis $\{\overline{s_n(x)}\}_{n\in{\mathbb N}}.$ It remains to note that, according to (\ref{0-a}), the knowledge
of $w_0(x)$ on $(\pi-a,\pi)$ is excessive since $w_1(x)$ has been found completely.

\medskip
{\bf Proof of Theorem~3.} Let us begin with the necessity part. According to (\ref{A2}), (\ref{CharInt-1}) and (\ref{CharInt-2}), we have
$$
\theta_0(\rho)=\int_0^\pi w_0(x)\sin\rho x\,dx, \quad \theta_1(\rho)=\int_0^\pi w_1(x)\cos\rho x\,dx.
$$
Hence, by virtue of (\ref{q+-}) and (\ref{0-a}), the exponential types of $\theta_0(\rho)$ and $\theta_1(\rho)$ do not exceed $\pi-a.$
Finally, the relation (\ref{A4}) follows from Lemmas~5 and~7 after substituting $\lambda=0$ into~(\ref{CharInt-2}) and the second formula in
(\ref{A3}). Indeed, according to (\ref{0-a}) and (\ref{a-pi-1}), the assumption $q^-=0$ implies
\begin{equation}\label{5.16}
\int_0^\pi w_1(x)\,dx =\frac14\int_a^\pi \Big(q^+\Big(\frac{a+x}2\Big) +q^+\Big(\pi+\frac{a-x}2\Big)\Big)\,dx =\frac12\int_a^\pi
q^+(x)\,dx=\omega.
\end{equation}

For the sufficiency, we construct the functions $\Delta_0(\lambda)$ and $\Delta_1(\lambda)$ by the formulae in~(\ref{A3}) using the given
sequences $\{\lambda_{n,0}\}_{n\ge0}$ and $\{\lambda_{n,1}\}_{n\ge0}.$  By virtue of Lemma~8, these functions have the forms
(\ref{CharInt-1}) and (\ref{CharInt-2}), respectively, with some $w_0(x),w_1(x)\in L_2(0,\pi),$ which, in turn, vanish a.e. on $(\pi-a,\pi)$
by the first condition along with the Paley--Wiener theorem \cite{Lev}.

By virtue of Lemma~10, there exists a unique solution $q^+(x)\in L_2(a,\pi)$ of the system (\ref{a-pi-0}) and (\ref{a-pi-1}) with $q^-(x)=0.$
As in (\ref{5.16}), we calculate
$$
\tilde\omega:=\int_0^{\pi-a} w_1(x)\,dx =\frac12\int_a^\pi q^+(x)\,dx
$$
and, hence,
\begin{equation}\label{5.18}
\Delta_1(0)=\omega +\tilde\omega.
\end{equation}
On the other hand, the second formula in (\ref{A3}) and condition (\ref{A4}) imply $\Delta_1(0)=2\omega,$ which along with (\ref{5.18}) gives
$\tilde\omega=\omega.$ Consider the problems $B_0(q)$ and $B_1(q)$ with the potential
$$
q(x)=\left\{\begin{array}{cl}0, &x\in(0,a),\\[3mm]
q^+(x), & x\in(a,\pi).\end{array}\right.
$$
According to Lemma~5, $\Delta_0(\lambda)$ and $\Delta_1(\lambda)$ are their characteristic functions, respectively. Hence,
$\{\lambda_{n,j}\}_{n\ge0}$ is the spectrum of $B_j(q)$ for $j=0,1.$ $\hfill\Box$
\\

{\bf Funding.} This research was supported by Russian Science Foundation, Grant No. 22-21-00509, https://rscf.ru/project/22-21-00509/
\\

{\bf Acknowledgement.} The authors are grateful to Maria Kuznetsova for reading the manu\-script and making valuable comments.



\end{document}